\documentclass[10pt]{article}

\usepackage{amssymb, amsmath}

\newtheorem{Thm}{Theorem}

\newtheorem{Cor}{Corollary}
\newtheorem{Prop}{Proposition}

\newtheorem{Fact}{Fact}

\newenvironment{thm}[3][]{\pagebreak[3] \medskip {\bf \noindent Theorem}
                               {#1} {#2}{\it #3}
                             }{}

\newtheorem{Rmk}{Remark}
\def\qed{\kern 6pt\hbox{\vrule\vbox to 6pt{\hrule width
                          6pt\vfil\hrule}\vrule}}

\newcommand{\R}{\mathbb{R}}
\newcommand{\C}{\mathbb{C}}
\renewcommand{\a}{\alpha}
\renewcommand{\d}{\delta}
\newcommand{\e}{\varepsilon}
\newcommand{\ga}{\gamma}

\def\e{\varepsilon}

\def\P{\mathcal P}
\newcommand{\la}{\lambda}
\def\f{\varphi}

\def\supp{{\rm supp}\,}
\def\relint{{\rm relint}}

\def\di{{\rm dim}\,}

\title{Concentration inequalities for $s$-concave measures of dilations of Borel sets and applications}
\author{ M. Fradelizi}
\begin{document}
\maketitle

\begin{center}
Universit\'e Paris-Est\\
Laboratoire d'Analyse et de Math\'ematiques
Appliqu\'ees (UMR 8050)\\
Universit\'e de Marne la Vall\'ee\\ 
77454 Marne la Vall\'ee Cedex 2, France \\
matthieu.fradelizi@univ-mlv.fr
\end{center}

\begin{abstract}
We prove a sharp inequality conjectured by Bobkov on the  measure of dilations of Borel sets in $\mathbb{R}^n$
by a $s$-concave probability. 
Our result gives a common generalization of  an inequality of Nazarov, Sodin and Volberg 
and a concentration inequality of Gu\'edon. Applying our inequality to the level sets of functions 
satisfying a Remez type inequality, we deduce, as it is classical, that these functions enjoy dimension free
distribution inequalities and Kahane-Khintchine type inequalities with positive and negative exponent, 
with respect to an arbitrary $s$-concave probability. \\
\end{abstract}



\newpage

\section{Introduction}

The main purpose of this paper is to establish a sharp inequality, conjectured by Bobkov in \cite{B3}, comparing the measure of a Borel set 
in $\R^n$ with a $s$-concave probability and the measure of its dilation. 
Among the $s$-concave probabilities are the log-concave ones ($s=0$) and thus the Gaussian ones, so that it is expected that they satisfy good concentration inequalities and large and small deviations inequalities. This is indeed the case and these inequalities as well as Kahane-Khintchine type inequalities with positive and negative exponent are deduced. By using a localization theorem in the form given by Fradelizi and Gu\'edon in \cite{FG}, we exactly determine among $s$-concave probabilities $\mu$ on $\R^n$ and among Borel sets $F$ in $\R^n$, with fixed measure $\mu(F)$, what is the smallest measure of the $t$-dilation of $F$ (with $t>1$). This infimum is reached for a one-dimensional measure which is $s$-affine (see the definition below) and $F=[-1,1]$. In other terms, it gives a uniform upper bound for the measure of the complement of the dilation of $F$ in terms of $t,s$ and $\mu(F)$.

The resulting inequality applies perfectly to sublevel sets of functions satisfying a  
Remez inequality, {\em i.e.} functions such that the $t$-dilation of any of their sublevel 
sets is contained in another of their sublevel set in a uniform way (see section 2.3 below). 
The main examples of such functions $f$ are the seminorms ($f(x)=\|x\|_K$, 
where $K$ is a centrally symmetric convex set in $\R^n$), 
the real polynomials in $n$-variables ($f(x)=P(x)=P(x_1, \dots, x_n)$, with $P\in\R[X_1, \dots, X_n]$)
and more generally  the seminorms of vector valued polynomials in $n$-variables 
($f(x)=\|\sum_{j=1}^N P_j(x)e_j\|_K$, with $P_1,\dots,P_N\in\R[X_1, \dots, X_n]$ and 
$e_1, \dots,e_N\in\R^n$). Other examples are given in section 3. For these functions we get 
an upper bound for the measures of their sublevel sets in terms of the measure of other sublevel sets. This enables to deduce that they satisfy large deviation inequalities and Kahane-Khintchine type inequalities with positive exponent.
But the main feature of the inequality obtained is that it may also be read backward. Thus it also implies
small deviation inequalities and Kahane-Khintchine type inequalities with negative exponent.\\

Before going in more detailed results and historical remarks, let us fix the notations.
Given subsets $A$, $B$ of the Euclidean space $\R^n$ and $\la \in\R,$ we set $A + B = \{x + y ;\  x \in A, y \in B\}$, 
$\la A = \{\la x ;\ x \in A\}$ and $A^c=\{x\in\R^n ;\ x\notin A\}$. For all $s \in (-\infty, 1],$ we say that a measure $\mu$
in $\R^n$ is $s$-concave if the inequality
$$
\mu(\la A + (1-\la) B ) \ge [\la \mu^s(A) + (1-\la)
\mu^s(B) ]^{1/s}
$$
holds for all compact subsets $ A, B \subset \R^n$ 
such that $\mu(A) \mu(B) > 0$   and  all $\la \in [0,1]$. 
The limit case is interpreted by continuity, thus the right hand side of this
inequality is equal to 
$\mu^{\la}(A) \, \mu^{1-\la}(B)$ for $s = 0$.
Notice that an $s$-concave measure is $t$-concave for all $t\le s$.
For a probability $\mu,$
$\supp (\mu)$ denotes its support.
For $\ga \in (-1, +\infty]$, a function  $f: \R^n \to \R^+$
 is $\ga$-concave if the inequality 
$$
f(\la x + (1-\la) y ) \ge [\la f^\ga(x) + (1-\la)
f^\ga (y) ]^{1/\ga}
$$
holds for all $x$ and $y$ such that  $f(x) f(y) > 0$ and all $\la \in [0,1]$,
where the limit cases $\ga=0$ and $\ga=+\infty$ are also interpreted by continuity, for example the $+\infty$-concave functions are constant. 
The link between the $s$-concave probabilities and the
$\ga$-concave functions is described in the work of Borell \cite{Bor2}.

\begin{thm}
    [\cite{Bor2}]
    {}
    {
    Let $\mu$ be a measure in $\R^n,$ let $G$ be the
    affine hull of the support of $\mu$, set $d =
    \di G$ and $m$ the Lebesgue measure on $G.$
    Then  for $s \le 1/d$, $\mu$ is $s$-concave if and only if
    $d\mu = \psi dm$, where 
    $0 \le \psi \in L_{loc}^{1}(\R^n, dm)$ and $\psi$ is $\ga$-concave 
    with $\ga = s/(1-sd)\in (-1/d,+\infty]$.
       }
\end{thm}

\medskip
\noindent 
According to this theorem, we say that a measure $\mu$ is
$s$-affine when its density $\psi$ satisfies  that
$\psi^{\ga}$  (or $\log \psi$ if $s=\ga=0$) is affine on its convex support
with $\ga=s/(1-sd)$.
In  \cite{Bor1},  Borell started the study of concentration properties of $s$-concave probabilities. 
He noticed  that for any centrally symmetric convex set $K$  the inclusion $K^c\supset \frac{2}{t+1} (tK)^c + \frac{t-1}{t+1}K$ holds true. 
From the definition of $s$-concavity he deduced that for every $s$-concave measure $\mu$ 
\begin{eqnarray}\label{Borell}
\mu(K^c)\ge \left(\frac{2}{t+1} \mu\big((tK)^c\big)^s + \frac{t-1}{t+1}\mu(K)^s\right)^{1/s}.
\end{eqnarray}
From this very easy but non-optimal concentration inequality, 
Borell showed that seminorms satisfy large deviation 
inequalities and Kahane-Khintchine type inequalities  with positive exponent. 
The same method was pushed forward in 1999 by 
Lata{\l}a \cite{L} to deduce a small ball probability for symmetric convex sets which allowed 
him to get a Kahane-Khintchine inequality until the geometric mean.

In 1991, Bourgain \cite{Bou} used the Knothe map \cite{K} to transport sublevel sets of polynomials. 
He deduced that, with respect to $1/n$-concave measure on $\R^n$ 
({\em i.e.} uniform measure on convex bodies), the real polynomials in $n$-variables 
satisfy some non-optimal distribution and Kahane-Khintchine type inequalities  with positive exponent. 
The same method was used by Bobkov in \cite{B2} and recently in \cite{B3} to generalize the result of Bourgain to $s$-concave measures 
and arbitrary functions, by using a "modulus of regularity" associated to the function. But  the concentration inequalities obtained in all these results using Knothe transport map  are not optimal.

In 1993, Lov\'asz and Simonovits \cite{LS} applied the localization method (using bisection arguments) 
to get the sharp inequality between the measure of a symmetric convex set $K$ and its dilation, for a log-concave probability $\mu$ 
\begin{eqnarray}\label{LoSi}
\mu\big((tK)^c\big)\le\mu(K^c)^{\frac{t+1}{2}}.
\end{eqnarray}
This improves inequality (\ref{Borell}) of Borell in the case $s=0$.
The method itself was further developped in 1995 by Kannan, Lov\'asz and Simonovits \cite{KLS} in a form more easily applicable. 
In 1999, Gu\'edon \cite{G} applied the localization method of \cite{LS} to generalize inequality (\ref{LoSi}) to the case of $s$-concave probabilities,
getting thus a full extension of inequality (\ref{Borell}). Gu\'edon proved that if  $\mu(tK)<1$ then
\begin{eqnarray}\label{Gue}
\mu(K^c)\ge \left(\frac{2}{t+1} \mu\big((tK)^c \big)^s + \frac{t-1}{t+1}\right)^{1/s}
\end{eqnarray}
and deduced from it the whole range of sharp inequalities (large and small deviations and Kahane-Khintchine) for symmetric convex sets. In 2000, Bobkov \cite{B1} used the localization in the form given in \cite{KLS} and the result of Lata{\l}a \cite{L} to sharpen the result of Bourgain on polynomials, with log-concave measures and proved that polynomials satisfy a  Kahane-Khintchine inequality until the geometric mean. 
In 2000 (published in 2002 \cite{NSV1}), Nazarov, Sodin and Volberg used the same bisection method to prove a "geometric Kannan-Lov\'asz-Simonovits lemma"  for log-concave measures. They generalized inequality (\ref{LoSi}) to arbitrary Borel set
\begin{eqnarray}\label{NaSoVo}
\mu(F_t^c)\le\mu(F^c)^{\frac{t+1}{2}},
\end{eqnarray}
where $F_t^c$ is the complement of $F_t$, the $t$-dilation of $F$, which is defined by
$$
F_t =\left\{x\in\R^n ;\ \hbox{there exists an interval $I\ni x $ such that}\  |I|<\frac{t+1}{2}|F\cap I|\right\} ,
$$
where $|\cdot|$ denotes the Lebesgue measure.
Notice that this definition of $t$-dilation is not the original definition of Nazarov, Sodin and Volberg \cite{NSV1}. In the later, 
they introduced an auxiliary compact convex set $K$ and used $t$ instead of $\frac{t+1}{2}$.  
The definition given above is the complement of their original one inside $K$. 
 The interest of our definition is that this auxiliary set becomes useless. 
If $F$ is open then its $t$-dilation  is open,  if $F$ is a  Borel set then its $t$-dilation is analytic, hence universally measurable.  
The $t$-dilation is an affine invariant, {\em i.e.} for any affine transform $A:\R^n\to\R^n$, we have 
$(AF)_t=A(F_t)$.
Notice that the definition of the $t$-dilation is one-dimensional in the sense that, if we denote by 
$\mathcal{D}$ the set of affine lines in $\R^n$, then
$$
F_t=\bigcup_{D\in\mathcal{D}}(F\cap D)_t .
$$
In \cite{NSV1}, Nazarov, Sodin and Volberg also noticed that $t$-dilation is well suited for sublevel sets of functions
satisfying a Remez type inequality and  deduced from the concentration inequality (\ref{NaSoVo})
that these functions satisfy  the whole range of sharp inequalities (large and small deviations and Kahane-Khintchine).
The preprint \cite{NSV1} had  a large diffusion and interested many people. For example, 
Carbery and Wright \cite{CW} and Alexander Brudnyi \cite{Br3}   directly applied the localization as presented in \cite{KLS} 
to deduce distributional inequalities and Kahane-Khintchine type inequalities for  the norm of vector valued polynomials in $n$-variables and 
functions with bounded Chebyshev degree, respectively.\\

Our main result is the following theorem which extends inequality (\ref{Gue}) of Gu\'edon to arbitrary Borel sets 
(since as we shall see in section 2, if $F$ is a centrally symmetric convex set $K$ then $F_t =tK$) and  inequality (\ref{NaSoVo}) of Nazarov, Sodin and Volberg to the whole range of $s$-concave probabilities. 
It establishes a conjecture of Bobkov \cite{B3} (who also proved in \cite{B3} a weaker inequality). 

\begin{Thm}\label{Main} 
Let $F$ be a Borel set in $\R^n$ and $t>1$. Let $s\in(-\infty, 1]$ and $\mu$ be a 
$s$-concave probability. Let 
$$
F_t =\{x\in\R^n ;\ \hbox{there exists an interval $I\ni x$ such that }\  |I|<\frac{t+1}{2}|F\cap I|\} .
$$
If $\mu(F_t)<1$ then 
\begin{eqnarray}\label{Mainineq}
\mu(F^c)\ge \left(\frac{2}{t+1} \mu(F_t^c)^s + \frac{t-1}{t+1}\right)^{1/s}.
\end{eqnarray}
\end{Thm}
Notice that inequality (\ref{Mainineq}) is sharp. For example, there is equality in (\ref{Mainineq})
if  $n=1$, $F=[-1,1]$ and $\mu$ is of density 
$$
\psi(x)=\frac{(a-s x)_+^{\frac{1}{s}-1}}{(a+s)^{\frac{1}{s}}}{\bf 1}_{[-1, +\infty)}(x) ,\ {\rm with}\ a>\max(-s,s t),
$$
with respect to the Lebesgue measure on $\R$, where $a_+=\max(a,0)$,  for every $a\in\R$.
Notice that this measure $\mu$ is $s$-affine on its support (which is $[-1, a/s]$ if $s>0$ 
and $[-1, +\infty)$ if $s\le 0$).\\

As noticed by Bobkov in \cite{B3}, in the case $s\le 0$, the right hand side term in inequality~(\ref{Mainineq}) 
vanishes if $\mu(F_t^c)=0$ so the condition $\mu(F_t)<1$ may be cancelled. 
But in the case $s>0$, the situation changes drastically. 
This condition is due to the fact that a $s$-concave probability measure, with $s>0$,
has necessarily a bounded support. From this condition we directly deduce the following corollary,
which was noticed by Gu\'edon \cite{G} in the case where $F$ is a centrally symmetric convex set.

\begin{Cor}\label{supp}
Let $F$ be a Borel set in $\R^n$. Let $s\in(0, 1]$ and $\mu$ be a 
$s$-concave probability. Denote by $V$ the relative interior of the (convex compact) support of $\mu$. Then
$$
V\subset F_t\quad \hbox{for every }\quad t\ge\frac{1+\mu(F^c)^s}{1-\mu(F^c)^s}.
$$
\end{Cor}

\noindent
{\bf Proof:} From Theorem \ref{Main}, if $\mu(F_t)<1$ then 
$$
\mu(F^c)\ge \left(\frac{2}{t+1} \mu(F_t^c)^s + \frac{t-1}{t+1}\right)^{1/s}> \left(\frac{t-1}{t+1}\right)^{1/s},
$$
which contradicts the hypothesis on $t$. Hence $\mu(F_t)=1$. It follows that $V\subset F_t$.\\

In section 2, we determine the effect of dilation on examples. The case of convex sets is treated in section 2.1, the case of 
sublevel sets of the seminorm of a vector valued polynomial in section 2.2 and the case of sublevel sets of a Borel measurable function
in section 2.3. In section 2.3, we also give a functional version of Theorem \ref{Main} and  we investigate the relationship between 
Remez inequality and inclusion of sublevel sets. In section 3 we deduce distribution and Kahane-Khintchine inequalities for functions 
of bounded Thebychev degree. Section 4 is devoted to the proof of Theorem \ref{Main}. The main tool for the proof is the localization 
theorem in the form given by Fradelizi and Gu\'edon in \cite{FG}.\\

After we had proven these results, we learned from Bobkov that, using a different method, Bobkov and Nazarov \cite{BN} simultaneously and independently proved Theorem \ref{Main}.\\

\section{Dilation of a set on examples}

\subsection{Convex sets}

\begin{Fact}   
Let $K$ be an open convex set then, for every $t>1$,
\begin{eqnarray}\label{eqconvex}
K_t=K+ \frac{t-1}{2}(K-K)=\frac{t+1}{2}K+\frac{t-1}{2}(-K)
\end{eqnarray}
and 
if moreover $K$ is centrally symmetric then $K_t=tK$.
\end{Fact}

\noindent
{\bf Proof:}
The second equality in (\ref{eqconvex}) deduces from  the convexity of $K$.
To prove the equality of the sets in (\ref{eqconvex}), we prove both inclusions:\\
Let $x\in K_t$. From the definition of $K_t$ and the remark following it, there exists a point $a\in\R^n$ such that
$|[a,x]|<\frac{t+1}{2}|K\cap[a,x]|$. Since $K$ is convex it follows that $K\cap[a,x]$ is an interval $[b,c]$. 
Let us denote the Euclidean norm by $|\cdot |_2$. Then $|x-a|_2<\frac{t+1}{2}|c-b|_2$, hence $b\neq c$. 
We may assume that $c\in (b,x]$ and $b\in[a,c]$. Hence there is 
$\lambda\in(0,1]$ such that $c=(1-\lambda)b+\lambda x$. This gives
$$
\frac{|c-b|_2}{\lambda}=|x-b|_2\le|x-a|_2<\frac{t+1}{2}|c-b|_2\ .
$$
Thus $\frac{1}{\lambda}<\frac{t+1}{2}$. Therefore
$$
x=c+\left(\frac{1}{\lambda}-1\right)(c-b)\in K+\left(\frac{1}{\lambda}-1\right)(K-K)\subset K+\frac{t-1}{2}(K-K)\ .
$$
Conversely, let $x\in \frac{t+1}{2}K+\frac{t-1}{2}(-K)$. If $x\in K$, the result is obvious so we assume that $x\notin K$.
There exists $b,c\in K$ such that $x=\frac{t+1}{2}c+\frac{t-1}{2}(-b)$.  Since $K$ is convex we deduce that the set 
$[b,x]\cap K$ is an interval with $b$ as an endpoint. 
Since $K$ is open there exists $d\in\R^n$ such that $[b,x]\cap K=[b,d)$ and we have $c\in [b,d)$.
Then 
$$
|[b,x]|=|x-b|_2=\frac{t+1}{2}|c-b|_2<\frac{t+1}{2}|d-b|_2=\frac{t+1}{2}\left|K\cap [b,x]\right|. 
$$
Therefore $x\in K_t$.\\
If  moreover $K$ is centrally symmetric it is obvious that $K_t=tK$. \hfill\qed\\

\noindent
{\bf Remarks:} \\
1) It is not difficult to see that if we only assume that $K$ is convex (and not necessarily open) then the same proof shows actually that  
$$K_t=\relint\left(K+ \frac{t-1}{2}(K-K)\right)=\relint\left(\frac{t+1}{2}K+\frac{t-1}{2}(-K)\right),$$ 
where $\relint(A)$ is the relative interior of $A$, {\em i.e.} the interior of $A$ relative to its affine hull.\\

\noindent
2)  The family of convex sets described by (\ref{eqconvex}) where introduced 
by Hammer \cite{H}, they may be equivalently defined in the following way. 
Let us recall that the support function of a convex set $K$ in the direction $u\in S^{n-1}$
is defined by $h_K(u)=\sup_{x\in K}\langle x,u\rangle$ and that an open convex set $K$ is equal to the 
intersection of the open slabs containing it:
$$
K=\bigcap_{u\in S^{n-1}}\!\left\{x\in\R^n;\ -h_K(-u)<\langle x,u\rangle < h_K(u)\right\}.
$$
The width of $K$ in direction $u\in S^{n-1}$ is defined by $w_K(u)=h_K(u)+h_K(-u)$.  
Then for every $t>1$,
$$
K_t=\!\!\bigcap_{u\in S^{n-1}}\left\{x;\ -h_K(-u)-\frac{t-1}{2}w_K(u)
<\langle x,u\rangle< h_K(u)+\frac{t-1}{2}w_K(u)\right\}.
$$
Moreover, since this definition can be extended to the values $t\in (0,1]$, it enables thus to define 
the $t$-dilation of a convex set for $0<t\le 1$ and in the symmetric case, the equality $K_t=tK$ is still valid for 
$t\in(0,1]$. Using that the family of convex sets $(K_t)_{t>0}$ is absorbing, Minkowski defined
what is now  called the "generalized Minkowski functional" of $K$:
 $$
\a_K(x)=\inf\left\{t>0;\ x\in K_t\right\}
$$
Notice that $\a_K$ is convex and positively homogeneous.
If moreover $K$ is centrally symmetric then $K_t=tK$, which gives $\a_K(x)=\|x\|_K$. 
We shall see in the next section, below,  how this notion was successfully used in 
polynomial approximation theory (see for example \cite{RS}).\\

From Fact 1, Theorem \ref{Main} and Corollary 1, we deduce the following corollary.

\begin{Cor}
Let $K$ be a convex set in $\R^n$ and $t>1$. Let $s\in(-\infty, 1]$ and $\mu$ be a 
$s$-concave probability. Let $V=\relint\big(\supp(\mu)\big)$.\\

\noindent
i) If $\mu\big(K+ \frac{t-1}{2}(K-K)\big)<1$ then 
$$
\mu(K^c)\ge \left(\frac{2}{t+1} \mu\Big(\big(K+ \frac{t-1}{2}(K-K)\big)^c \Big)^s + \frac{t-1}{t+1}\right)^{1/s}.
$$
ii) If $s>0$ then 
$$
V\subset K+\frac{\mu(K^c)^s}{1-\mu(K^c)^s}(K-K)
$$
iii) If $s>0$ and $K$ is centrally symmetric then
$$
V\subset\frac{1+\mu(K^c)^s}{1-\mu(K^c)^s}K
$$ 
\end{Cor}

\noindent
Applying $iii)$ to the uniform probability on $V$ we deduce  that for every convex sets $V$ and $K$ in $\R^n$, with $K$ symmetric
$$
V\subset\frac{|V|^\frac{1}{n}+|V\cap K^c|^\frac{1}{n}}{{|V|^\frac{1}{n}-|V\cap K^c|^\frac{1}{n}}}K.
$$

\subsection{Sublevel set of the seminorm of a vector valued polynomial}

Let $P$ be a polynomial of degree $d$, with $n$ variables and with values in a Banach space $E$, that is 
$$
P(x_1,...,x_n)=\sum_{k=1}^NP_k(x_1,...,x_n)e_k\ ,
$$
where $e_1, ...,e_N\in E$ and $P_1,...,P_N$ are real polynomials with $n$ variables and degree at most $d$.
Let  $K$ be a centrally symmetric convex set in $E$, and denote by $\|\cdot\|_K$ the seminorm defined by $K$ in $E$ and let $c>0$ be any constant. The following fact  was noticed and used 
by Nazarov, Sodin and Volberg in \cite{NSV1}, in the case of real polynomials.

\begin{Fact}
Let $P$ be a polynomial of degree $d$, with $n$ variables and with values in a Banach space $E$
 and let $t>1$. Let  $K$ be a centrally symmetric convex set in $E$ and $c>0$. Then
$$
\{x\in\R^n;\ \|P(x)\|_K<c\}_t\subset\{x\in\R^n;\ \|P(x)\|_K<c\,T_d (t)\},
$$
where $T_d$ is the Chebyshev polynomial of degree $d$, {\em i.e.}
$$
T_d(t)=\frac{\left(t+\sqrt{t^2-1}\right)^d+\left(t-\sqrt{t^2-1}\right)^d}{2},
$$
for every $t\in\R$ such that $|t|\ge 1$. 
\end{Fact}

This fact is actually a reformulation, in terms of dilation, of the Remez inequality \cite{R} which asserts 
that for every real polynomial $Q$ of degree $d$ and one variable, for every interval $I$ in $\R$ and every Borel subset 
$J$ of $I$,  
$$
\sup_I |Q|\le T_d\left(2\frac{|I|}{|J|}-1\right)\sup_J|Q| .
$$
Let us prove the inclusion. Let $x_0\in F_t$. There exists an interval $I=[a,b]$ containing $x_0$ such that $|I|<\frac{t+1}{2}|F\cap I|$.
The key point is that 
$$
\|P\big((1-\lambda)a+\lambda b\big)\|_K=\sup_{\xi\in K^*}\xi\Big(P\big((1-\lambda)a+\lambda b\big)\Big)=\sup_{\xi\in K^*}Q_\xi(\lambda),
$$
where $K^*=\{\xi\in E^*; \forall\ x\in K, \xi(x)\le 1\}$ is the polar of $K$ and 
$Q_\xi(\lambda)= \xi\Big(P\big((1-\lambda)a+\lambda b\big)\Big)$ is a real polynomial of one 
variable and degree at most $d$. 
Let $J:=\{\lambda\in[0,1];\ (1-\lambda)a+\lambda b\in F\}$, then $|J|=|F\cap I|/|I|$.
Applying Remez inequality to $Q_\xi$ we have
$$
\sup_{\lambda\in[0,1]}Q_\xi(\lambda)\le T_d\left(\frac{2}{|J|}-1\right)\sup_{\lambda\in J}|Q_\xi(\lambda)| 
=T_d\left(\frac{2|I|}{|F\cap I|}-1\right)\sup_{x\in F\cap I}|\xi\big(P(x)\big)| .
$$
Taking the supremum, using that $T_d$ is increasing on $[1, +\infty)$ and the definition of $F$,  we get 
\begin{eqnarray*}
\|P(x_0)\|_K & \le & \sup_{[0,1]}\|P\big((1-\lambda)a+\lambda b\big)\|_K=\sup_{[0,1]}\sup_{\xi\in K^*}Q_\xi(\lambda)\\
& \le & T_d\left(\frac{2|I|}{|F\cap I|}-1\right)\sup_{\xi\in K^*}\sup_{x\in F\cap I}\xi\big(P(x)\big) \\
& < & T_d(t) \sup_{x\in F\cap I}\|P(x)\|_K\le cT_d(t) .
\end{eqnarray*}

\noindent
{\bf Remark:}
Notice that the Chebyshev polynomial of degree one is $T_1(t)=t$. Hence if we take the polynomial 
$P(x)=x=\sum x_i e_i$, where $(e_1, ..., e_n)$ is the canonical orthonormal basis of $\R^n$, we see 
that the case of vector valued polynomials generalizes the case of symmetric convex sets. \\

Fact 2 has an interesting reformulation in terms of 
polynomial inequalities in real approximation theory. It may be written in the following way.
Denote by $\mathcal{P}_d^n(E)$ the set of  polynomials of degree $d$, with $n$ variables 
and with values in a Banach space $E$.
Let $P\in\mathcal{P}_d^n(E)$ and $K$ be a symmetric convex set in $E$. 
Let $F$ be a Borel set in $\R^n$ and $t> 1$. For $x\in F_t$
$$
\|P(x)\|_K\le T_d(t)\sup_{z\in F} \|P(z)\|_K .
$$
Let us assume that the Borel set $F$ in $\R^n$ has the property that, 
for each $x$ in $R^n$, there is an affine line $D$ containing $x$ 
such that $|F\cap D|>0$, which is the case if $F$ has non-empty interior.
Then $\bigcup_{t>1}F_t=\R^n$.
In this case, we may define for every $x\in\R^n$ the
"generalized Minkowski functional" of $F$ at $x$ as
$$
\a_F(x)=\inf\{t>1;\ x\in F_t\}.
$$
Using this quantity, we get the following reformulation of Fact 2.

\begin{Cor}\label{corfact2}
Let $F$ be a Borel set in $\R^n$. Let $P\in\mathcal{P}_d^n(E)$ and 
$K$ be a centrally symmetric convex set in $E$. For every $x$ in $\R^n$, 
$$
\|P(x)\|_K\le T_d\big(\a_F(x)\big)\sup_{z\in F} \|P(z)\|_K .
$$
\end{Cor}

Let us introduce the notations coming from approximation theory. With the notations of the 
corollary, we define 
$$
C_d(F,x,K)=\sup\{\|P(x)\|_K; P\in\mathcal{P}_d^n(E), \sup_{x\in F}\|P(x)\|_K\le 1, n\ge 1 \}.
$$
Then the inequality may be written in the following form.
$$
C_d(F,x,K)=T_d\big(\a_F(x)\big).
$$
For $F$ being convex and the polynomial $P$ being real valued, this is a theorem of 
Rivlin-Shapiro \cite{RS} (see also an extension in \cite{RSa1} and  \cite{RSa2}). We get thus an extension
of their theorem to non-convex sets $F$, as well as Remez inequality generalizes to Borel sets
the classical Tchebychef inequality valid for segments in $\R$. \\

Applying Theorem \ref{Main} to the level set of a polynomial we get the following corollary, which was proved in the case $s=0$ by Nazarov, Sodin and Volberg in \cite{NSV1} and in the case $d=1$ and $P(x)=x$ by Gu\'edon in \cite{G}.

\begin{Cor}\label{mupoly}
Let $P$ be a polynomial of degree $d$, with $n$ variables and with values in a Banach space $E$
 and let $t>1$. Let $K$ be a centrally symmetric convex set in $E$ and $c>0$. Let $s\le1$ and $\mu$ be a $s$-concave probability. 
 If $\mu(\{x;\ \|P(x)\|_K\ge cT_d(t)\})>0$, then
$$
\mu(\{x;\ \|P(x)\|_K\ge c\})\ge\left(\frac{2}{t+1} \mu(\{x;\ \|P(x)\|_K\ge cT_d(t)\})^s + \frac{t-1}{t+1}\right)^{1/s}.
$$
For $s=0$, 
$$
 \mu(\{x;\ \|P(x)\|_K\ge cT_d(t)\})\le \mu(\{x;\ \|P(x)\|_K\ge c\})^{\frac{t+1}{2}}.
$$
\end{Cor}

Applying Corollary 1, we get the following extension of a theorem of Brudnyi and Ganzburg \cite{BG} (which treats the case of probabilities $\mu$ which are  uniform on a convex body). It is a multi-dimensional version of Remez inequality.

\begin{Cor}\label{suppoly}
Let $P$ be a polynomial of degree $d$, with $n$ variables and with values in a Banach space $E$. Let $K$ be a centrally symmetric convex set in $E$. 
Let $s\in(0,1]$, $\mu$ be a $s$-concave probability and let $V$ be the support of $\mu$.  Then, for every $\omega\subset V$
$$
\sup_{x\in V}\|P(x)\|_K\le T_d\left(\frac{1+\mu(\omega^c)^s}{1-\mu(\omega^c)^s}\right)\sup_{x\in\omega}\|P(x)\|_K
\le\left(\frac{4}{s\mu(\omega)}\right)^d\sup_{x\in\omega}\|P(x)\|_K.
$$
\end{Cor}
 
\noindent
{\bf Proof:}
We apply Corollary 1 to $F=\{x ;\ \|P(x)\|_K\le \sup_{x\in\omega}\|P(x)\|_K\}$ and Fact 2 to deduce that 
$$
V\subset F_t\subset \{x\in\R^n;\ \|P(x)\|_K<T_d (t)\sup_{x\in\omega}\|P(x)\|_K\},\quad \forall\ t\ge\frac{1+\mu(F^c)^s}{1-\mu(F^c)^s}.
$$
Since $\omega\subset F$, we may apply the preceding inclusion to $t=\frac{1+\mu(\omega^c)^s}{1-\mu(\omega^c)^s}$ and this gives the first inequality.
The second one follows using that $T_d(t)\le (2t)^d$ for every $t\ge 1$ and easy computations.

\subsection{Sublevel set of a Borel measurable function}

In Fact 1 and Fact 2 we saw the effect of dilation on convex sets and level sets of vector valued polynomials. We want to describe now the most general case of level sets of Borel measurable functions. As in Fact 2, we shall see in the following proposition that for any Borel measurable function,
 an inclusion between the dilation of the level sets is equivalent to a Remez type inequality.

 \begin{Prop}\label{propuf}
Let  $f:\ \R^n\to \R$ be a Borel measurable function and $t>1$. 
Let $u_f(t)\in [1,+\infty)$. The following are equivalent.\\
i) For every interval $I$ in $\R^n$ and every Borel subset $J$ of $I$ such that $|I|<t|J|$,  
$$
\sup_I|f|\le u_f(t)\sup_J|f| .
$$
ii) For every $\la>0$, 
$$
\left\{x\in\R^n;\ |f(x)|\le \la\right\}_{2t-1}\subset \left\{x\in\R^n;\ |f(x)|\le \la u_f(t)\right\}.
$$
\end{Prop}

We shall say that a non-decreasing function $u_f:\ (1,+\infty)\to [1,+\infty)$ is {\it a Remez function} of $f$ 
if it satisfies i) or ii) of the previous proposition, for every $t>1$ and that it is  {\it the Remez function} of $f$ if it is the smallest Remez function of $f$.  

For example, using {\it i)}, the Remez inequality asserts that if we take $f(x)=\|P(x)\|_K$ 
where $P$ is a polynomial of degree $d$, with $n$ variables and 
with values in a Banach space $E$ and $K$ is a symmetric convex set 
then $t\mapsto T_d(2t-1)$ is a Remez function of $f$. Using  {\it ii)} and Fact 1,  we get that 
$u_f(t)=2t-1$ is the Remez function of $f(x)=\|x\|_K$.\\

\noindent
{\bf Proof of Proposition \ref{propuf}:}\\
$i)\Longrightarrow ii)$:
Let $F=\{x\in\R^n;\ |f(x)|\le\la\}$ and let $x\in F_{2t-1}$. 
There exists an interval $I$ containing $x$ such that $|I|<t|F\cap I|$. Hence 
$$
|f(x)|\le\sup_I|f|\le u_f(t)\sup_{F\cap I}|f|\le \la u_f(t).
$$
$ii)\Longrightarrow i)$:
Let $I$ be an interval in $\R^n$ and $J$ be a Borel subset of $I$
such that  $|I|< t|J|$. Let $\la=\sup_J|f|$ and let $x\in I$, then $J\subset \{|f|\le \la\}\cap I $ hence 
$$
 |I|<t|J|\le t|\{|f|\le \la\}\cap I|,
$$
thus $x\in \{|f|\le \la\}_{2t-1}$. From $ii)$ we get $|f(x)|\le \la u_f(t)$. This gives $i)$.
\hfill\qed\\

\noindent
Applying Theorem \ref{Main} to the level set of a Borel measurable function, we get the following.

\begin{Thm}\label{Mainfunc} 
Let  $f:\R^n\to \R$ be a Borel measurable function and
 $u_f:\ (1,+\infty)\to [1,+\infty)$ be a Remez function of $f$. 
Let $s\in(-\infty, 1]$ and $\mu$ be a $s$-concave probability. 
Let $t>1$ and $\la>0$. If $\mu(\{x;\ |f(x)|\ge \la u_f(t)\})>0$, then
\begin{eqnarray}\label{Mainfuncineq}
\mu(\{x;\ |f(x)|> \la\})\ge\left(\frac{1}{t} \mu(\{x;\ |f(x)|> \la u_f(t)\})^s + 1-\frac{1}{t}\right)^{1/s}.
\end{eqnarray}
For $s=0$,
$$
\mu(\{x;\ |f(x)|> \la u_f(t)\})\le \mu\left(\{x;\ |f(x)|> \la\}\right)^t.
$$
\end{Thm}

\noindent
{\bf Remark:} 
Theorem \ref{Mainfunc} improves a theorem given by Bobkov in \cite{B3}. 
As in \cite{B3}, notice that Theorem \ref{Mainfunc} is a functional version of Theorem \ref{Main}. 
As a matter of fact, we may follow the proof given by Bobkov. 
If a Borel subset $F$ of $\R^n$ and $u>1$ are given,
we apply Theorem~\ref{Mainfunc} to $t=\frac{u+1}{2}$,  $\la=1$ and 
$$
f=1 \quad \hbox{on}\ F,\quad f=2 \quad \hbox{on}\ F_u\setminus F\quad\hbox{and}\quad f=4\quad \hbox{on}\  F_u^c.
$$
Using ii) of Proposition \ref{propuf} it is not difficult to see that $u_f(t)=2$. Then 
inequality~(\ref{Mainineq}) follows from inequality~(\ref{Mainfuncineq}).\\

Applying Corollary 1 in the similar way as in Corollary \ref{suppoly} and using Proposition \ref{propuf} instead of Fact 2,  we get the following.

\begin{Cor}
Let  $f:\R^n\to \R$ be a Borel measurable function. 
Let $u_f:\ (1,+\infty)\to [1,+\infty)$ be a Remez function of $f$. 
Let $s\in(0, 1]$ and $\mu$ be a $s$-concave probability. 
Let $\omega\subset\R^n$, then
$$
\| f\|_{L^\infty(\mu)}\le\sup_\omega |f|\ u_f\left(\frac{1}{1-\mu(\omega^c)^s}\right)\le\sup_\omega |f|\ u_f\left(\frac{1}{s\mu(\omega)}\right).
$$
\end{Cor}

Instead of using $u_f$,  Bobkov in \cite{B2} and \cite{B3}  introduced  a related quantity,
the "modulus of regularity" of $f$,  
$$
\d_f(\e)=\sup_{x,y}|\{\la\in[0,1];\ |f((1-\la)x+\la y)|\le\e|f(x)|\}|, \quad {\rm for}\ 0<\e\le1.
$$
It is not difficult to see that 
$$
\d_f(\e)=\sup_{x,y}\frac{|\{z\in[x,y];\ |f(z)|\le\e\sup_{[x,y]}|f|\}|}{|[x,y]|}
$$
and thus
$$
\d_f(\e)=\sup\left\{\frac{|J|}{|I|};\ J\subset I, \hbox{where $I$ is an interval and} 
\sup_J|f|\le\e\sup_I|f| \right\}.
$$
Hence $\d_f$ is the smallest function satisfying that for every interval $I$ and 
every Borel subset $J$ of $I$
$$
\frac{ |J|}{|I|}\le\d_f\left(\frac{ \sup_J |f|}{\sup_I|f|}\right),
$$
which is a  Remez-type inequality.
For smooth enough functions, the relationship between $u_f$, the Remez function of $f$ and $\d_f$ is given by
$$
\d_f(\e)=\frac{1}{u_f^{-1}(1/\e)},
$$
where $u_f^{-1}$ is the reciprocal function of $u_f$.
Hence if $f(x)=\|P(x)\|_K$ where $P$ is a polynomial 
of degree $d$, with $n$ variables and 
with values in a Banach space $E$ and $K$ is a symmetric convex set 
then, using that $u_f(t)\le T_d(2t-1)$ and $T_d(t)\le 2^{d-1}t^d$, for every $|t|\ge 1$, we get
$$
u_f(t)\le T_d(2t-1)\le 2^{d-1}(2t-1)^d\le (4t)^d
$$
and
\begin{eqnarray}\label{eqndelta}
\d_f(\e)\le\frac{2}{T_d^{-1}(1/\e)+1}\le 4\left(\frac{\e}{2}\right)^{1/d}\le 4\e^{1/d},
\end{eqnarray}
 for every $|t|\ge 1$. For $f(x)=\|x\|_K$, we get 
$\d_f(\e)=\frac{2\e}{\e+1}$ as noticed by Bobkov in \cite{B2}.
Notice that inequalities (\ref{eqndelta})  improve the previous bound given by
Bobkov in \cite{B2} and \cite {B3}. 

The interest of the quantity $\d_f$ comes from the next corollary, which 
 was conjectured by Bobkov in \cite{B3} 
 (for $s=0$, it deduces from \cite{NSV1}  as noticed in \cite{B2}).

\begin{Cor}\label{mufonc}
Let  $f:\ \R^n\to \R$ be a Borel measurable function and $0<\e\le 1$. 
Let $s\le1$ and $\mu$ be a $s$-concave probability. Let $\la<\|f\|_{L^\infty(\mu)}$, then
\begin{eqnarray}\label{ineqdelta}
\mu(\{|f|\ge \la\e\})\ge\left(\d_f(\e) \mu(\{|f|\ge \la\})^s + 1-\d_f(\e)\right)^{1/s}
\end{eqnarray}
and if $\mu$ is $\log$-concave ({\em i.e.} for $s=0$) then
$$
\mu(\{|f|\ge \la\e\})\ge \mu(\{|f|\ge \la\})^{\d_f(\e)} .
$$
\end{Cor}

{\bf Proof:}
We apply Theorem \ref{Main} to the set $F=\{|f|<\la\e\}$ and $t=\frac{2}{\d_f(\e)}-1$.
Let $x\in F_t$, there exists an interval $I$ containing $x$ such that 
$$
|I|<\frac{t+1}{2}|F\cap I|=\frac{|F\cap I|}{\d_f(\e)}.
$$
From the definition of $\d_f$, this implies that  
$$
f(x)\le\sup_I|f|<\frac{1}{\e}\sup_{F\cap I}|f|\le\la.
$$
Hence $F_t\subset\{|f|< \la\}$. This gives the result.\hfill\qed

\section{Distribution and Kahane-Khintchine type inequalities}

It is classical that from an inequality like inequality (\ref{Mainfuncineq}) (or in its equivalent form (\ref{ineqdelta})), it is possible to deduce distribution and Kahane-Khintchine type inequalities. Due to its particular form, this type of concentration inequality may be read forward or backward and thus permits  to deduce both small and large deviations inequalities. 

\subsection{Functions with bounded Chebyshev degree}

Before stating these inequalities, let us define an interesting  set of functions, 
the functions $f$ such that their Remez function $u_f$ is bounded from above by a power function, {\em i.e.} 
 there exists $A>0$ and $d>0$ satisfying  $u_f(t)\le (At)^{d}$, for every $t>1$
which means that for every interval $I$ in $\R^n$ and every Borel subset $J$ of $I$
$$
\sup_I|f|\le \left(\frac{A|I|}{|J|}\right)^{d}\sup_J|f| .
$$
In this case, the smallest power satisfying this inequality is called {\it the Chebyshev degree} of $f$ and denoted by $d_f$ and the best constant 
corresponding to this degree is denoted by $A_f$. 
This is also equivalent to assume that $\d_f(\e)\le A_f\e^{1/d_f}$, for every $0<\e<1$.
Notice that if $f$ has bounded Chebyshev degree ({\it i.e.} $d_f <+\infty$)  then
$|f|^{1/d_f}$ has Chebyshev degree one and $A_{|f|^{1/d_f}}=A_f$. For such functions inequality (\ref{Mainfuncineq})
becomes, for every $t>1$,
\begin{eqnarray}\label{ineqTcheb}
\mu(\{|f|^{1/d_f}>\la\})\ge\left(\frac{1}{t} \mu(\{x;\ |f(x)|^{1/d_f}> \la A_f t\})^s + 1-\frac{1}{t}\right)^{1/s}
\end{eqnarray}
and for $s=0$
\begin{eqnarray}\label{logTcheb}
\mu(\{x;\ |f(x)|^{1/d_f}> \la A_f t\})\le\mu(\{|f|^{1/d_f}>\la\})^t.
\end{eqnarray}
 For example if $f(x)=\|x\|_K$ then $u_f(t)=2t-1$ hence $d_f=1$ and $A_f=2$. 
 If $f(x)=\|P(x)\|_K$ where $P$ is a polynomial of degree $d$, with $n$ variables and 
with values in a Banach space $E$ and $K$ is a symmetric convex set  then  $d_f=d$ and $A_f=4$.
More generally, following \cite{NSV1} and \cite{CW}, if  $f=e^{u}$, where $u:\R^n\to\R$ is the restriction to $\R^n$ of a plurisubharmonic function $\tilde{u}:\C^n\to\R$ such that $\limsup_{|z|\to +\infty} \frac{\tilde{u}(z)}{\log |z|}\le 1$, then $d_f=1$ and $A_f=4$. 
Another type of example was given by Nazarov, Sodin and Volberg in \cite{NSV1}: if 
$$
f(x)=\sum_{k=1}^dc_k e^{i\langle x_k,x\rangle},
$$
with $c_k\in\C$ and $x_k\in\R^n$ then $d_f=d$. Finally, Alexander Brudnyi in \cite{Br1}, \cite{Br2}, \cite{Br3} 
(see also Nazarov, Sodin and Volberg \cite{NSV2}) proved also that for any $r>1$, 
for any holomorphic  function $f$ on $B_\C(0,r)\subset\C^n$, the open complex Euclidean of radius $r$ centered at $0$,
the Chebyshev degree of $f$ is bounded.

 \subsection{Small deviations and Kahane-Khintchine type inequalities for negative exponent}

Let us start with the following small deviation inequality, which was proved by Gu\'edon 
\cite{G} in the case where $f=\|\cdot\|_K$ and by  Nazarov, Sodin and Volberg \cite{NSV1}
in the case where $s=0$. It was  proved in a weaker form and conjectured in this form
by Bobkov in \cite{B3}.
This type of inequality is connected to small ball probabilities.

\begin{Cor}\label{smalldelta}
Let  $f:\ \R^n\to \R$ be a Borel measurable function and $0<\e\le 1$. 
Let $s\le1$ and $\mu$ be a $s$-concave probability. Let $\la<\|f\|_{L^\infty(\mu)}$, then
\begin{eqnarray}\label{ineqsmalldelta}
\mu(\{|f|\le \la\e\})\le\d_f(\e)\times\frac{1- \mu(\{|f|\ge \la\})^s }{s}.
\end{eqnarray}
In particular, if $\mu$ is $\log$-concave ({\em i.e.} for $s=0$) then
$$
\mu(\{|f|\le \la\e\})\le \d_f(\e)\times\log\left(1/ \mu(\{|f|\ge \la\})\right) .
$$
\end{Cor}

{\bf Proof:}
Let $s\neq 0$.The proof given by Gu\'edon in \cite{G} works here also. We reproduce it here for completeness.
Since $s\le1$ the function $x\mapsto (1-x)^{1/s}$ is convex on $(-\infty,1]$, hence  
$$
(1-x)^{1/s}\ge 1-\frac{x}{s}.
$$
The result follows from inequality (\ref{ineqdelta}) and the  inequality above applied to 
$x=\d_f(\e)(1-\mu(\{|f|\ge \la\})^s)$. For $s=0$ the result follows by taking limits or can be proved along the same lines.
\hfill\qed\\

In the case of functions with bounded Chebyshev degree, inequality (\ref{ineqsmalldelta}) take a simpler form  
and, by integrating on level sets, it immediately gives an inverse H\"older Kahane-Khintchine type inequality for negative exponent.
Thus, we get the following corollary, generalizing a theorem of Gu\'edon \cite{G} (for $f=\|.\|_K$) and Nazarov, Sodin and Volberg \cite{NSV1} (for $s =0$).

\begin{Cor}\label{Khintchnega}
Let  $f:\ \R^n\to \R$ be a Borel measurable function with bounded Chebyshev degree. 
Let $s\le1$ and $\mu$ be a $s$-concave probability. Denote by $M_f$ the $\mu$-median of $|f|^{1/d_f}$ 
and denote $c_s:=(1-2^{-s})/s$, for $s\neq 0$ and $c_0=\ln 2$. Then for every $0<\e <1$
\begin{eqnarray}\label{ineqsmalldf}
\mu(\{|f|^\frac{1}{d_f}\le M_f\e\})\le A_fc_s\e,
\end{eqnarray}
and for every $-1<q<0$, 
\begin{eqnarray}\label{Khintchneg}
\||f|^\frac{1}{d_f}\|_{L^q(\mu)}\ge M_f\left(1-\frac{qA_fc_s}{q+1}\right)^{1/q}\ge M_f e^{-\frac{A_fc_s}{q+1}}\ .
\end{eqnarray}
\end{Cor}

{\bf Proof:} Inequality (\ref{ineqsmalldf}) deduces from inequality (\ref{ineqsmalldelta}) by taking $\la=M_f$. 
The proof of inequality (\ref{Khintchneg}) is then standard, we apply  inequality (\ref{ineqsmalldf})
\begin{eqnarray*}
\int_{\R^n}|f(x)|^\frac{q}{d_f}d\mu(x) &=& -q\int_0^{+\infty}t^{q-1}\mu\left(\{x;\ |f(x)|^\frac{1}{d_f}\le t\}\right) dt\\
&\le& -q\int_0^{M_f}t^{q}\frac{A_fc_s}{M_f}dt -q\int_{M_f}^{+\infty}t^{q-1}dt\\
&=& M_f^q\left(1-\frac{qA_fc_s}{q+1}\right).
\end{eqnarray*}
Then we take the $q$-th root (recall that $q<0$) to get inequality (\ref{Khintchneg}).

 \subsection{Large deviations and Kahane-Khintchine type inequalities for positive exponent}

On the contrary to the small deviations case, the behaviour of large deviations of a function 
with bounded Chebyshev degree with respect to a $s$-concave probability heavily depends 
on the range of $s$, mainly on the sign of $s$.  But all behaviours follow from inequality  (\ref{ineqTcheb}) applied to $\la=M_f$, 
the $\mu$-median of $|f|^{1/d_f}$, which gives, for every $s\le 1$, $s\neq 0$,
\begin{eqnarray}\label{ineqlarge}
\mu(\{|f|^\frac{1}{d_f}\ge A_f M_ft\})\le\left(1-t(1-2^{-s})\right)_{+}^{\frac{1}{s}}
\end{eqnarray}
and for $s=0$,
$$
\mu(\{|f|^\frac{1}{d_f}\ge A_f M_ft\})\le 2^{-t}.
$$
For $s\ge 0$, it follows from inequality (\ref{ineqlarge}) that $|f|^{1/d_f}$ has exponentially decreasing tails
and a standard argument implies an inverse H\"older inequality.

\begin{Cor}\label{Khintchnposspos}
Let  $f:\ \R^n\to \R$ be a Borel measurable function with bounded Chebyshev degree. 
Let $0\le s\le1$ and $\mu$ be a $s$-concave probability. Denote by $M_f$ the $\mu$-median of $|f|^{1/d_f}$ 
and denote $c_s:=(1-2^{-s})/s$, for $s> 0$ and $c_0=\ln 2$. Then for every $t>1$
\begin{eqnarray}\label{ineqspos}
\mu(\{|f|^\frac{1}{d_f}\ge A_f M_ft\})\le\left(1-sc_st\right)_{+}^{\frac{1}{s}}\le e^{-c_s t}\le e^{-\frac{t}{2}}
\end{eqnarray}
and for every $p>0$, 
\begin{eqnarray}\label{ineqKhintposspos}
\||f|^\frac{1}{d_f}\|_{L^p(\mu)}\le A_f M_f\left(1+\frac{pB(p,1+\frac{1}{s})}{(sc_s)^p}\right)^{\frac{1}{p}}\le A_f M_f\left(1+2^p\Gamma(p+1)\right)^{\frac{1}{p}}.
\end{eqnarray}
\end{Cor}

{\bf Proof:} Inequality (\ref{ineqspos}) deduces from inequality (\ref{ineqlarge}). The proof of inequality (\ref{ineqKhintposspos}) is then standard, we write 
$$
\int_{\R^n}|f(x)|^\frac{p}{d_f}d\mu(x)=p\int_0^{+\infty}t^{p-1}\mu\left(\{x;\ |f(x)|^\frac{1}{d_f}\ge t\}\right) dt
$$
and we apply inequality (\ref{ineqspos}) as in the proof of Corollary \ref{Khintchnega}.\hfill\qed\\

For $s<0$ the situation changes drastically,  inequality (\ref{ineqlarge}) only implies that the tail of $|f|^{1/d_f}$ decreases as $t^{1/s}$, which is the sharp behaviour if we take the example of measure $\mu$ on $\R$ given after Theorem 1 and $f(x)=|x|$.

\begin{Cor}\label{Khintchnpossneg}
Let  $f:\ \R^n\to \R$ be a Borel measurable function with bounded Chebyshev degree. 
Let $s\le 0$ and $\mu$ be a $s$-concave probability. Denote by $M_f$ the $\mu$-median of $|f|^{1/d_f}$ and denote $d_s:=(2^{-s}-1)^{1/s}$. Then for every $t>1$
\begin{eqnarray}\label{ineqsneg}
\mu(\{|f|^{1/d_f}\ge A_f M_ft\})\le t^{\frac{1}{s}}\left(2^{-s}-1+\frac{1}{t}\right)^{\frac{1}{s}}\le d_s t^{\frac{1}{s}}.
\end{eqnarray}
and for every $0<p<-\frac{1}{s}$, 
\begin{eqnarray}\label{ineqKhintnegspos}
\||f|^{1/d_f}\|_{L^p(\mu)}\le A_f M_f\left(1+d_s \frac{p }{p+\frac{1}{s}}\right)^{\frac{1}{p}}.
\end{eqnarray}
\end{Cor}

{\bf Proof:} Inequality (\ref{ineqsneg}) deduces from inequality (\ref{ineqlarge}). The proof of inequality (\ref{ineqKhintnegspos}) is then standard.\\

\section{Proof of Theorem \ref{Main}}

While in \cite{B2} and \cite{B3}, Bobkov used an argument based on a transportation argument, 
going back to Knothe \cite{K} and Bourgain \cite{Bou},
our proof follows the same line of argument as Lov\'asz and Simonovits in \cite{LS}, Gu\'edon in \cite{G}, Nazarov, Sodin and Volberg in \cite{NSV1}, Brudnyi in \cite{Br3} and Carbery and Wright in \cite{CW}, 
the geometric localization theorem, which reduces the problem to the dimension one. The main difference with these 
 proofs is that the geometric localization is used here in the presentation given by Fradelizi and Gu\'edon in \cite{FG} which don't use an infinite bisection method but prefers to see it as an optimization theorem on the set of $s$-concave measures satisfying a linear constraint and the application of the Krein-Milman theorem.
Let us recall the main theorem of \cite{FG}.

\begin{thm}
    [\cite{FG}] {}{Let $n$ be a positive integer, let $K$ be a compact convex set in $\R^n$ and denote by $\P(K)$ the set of probabilities in $\R^n$ supported in $K$. Let $f:K\to\R$ be an upper semi-continuous function,
let $s\in[-\infty, \frac{1}{2}]$ and denote by $P_f$ the set of $s$-concave probabilities $\la$ supported in $K$ satisfying $\int fd\la\ge 0$. Let $\Phi: \P(K)\to\R$ be a convex upper semi-continuous function. Then
$$
\sup_{\la\in P_f}\Phi(\la)
$$
is achieved at a probability $\nu$ which is either a Dirac measure at a point $x$ such that $f(x)\ge0$, 
or a probability $\nu$ which is $s$-affine on a segment $[a,b]$, such that $\int fd\nu=0$ and 
$\int_{[a,x]}fd\nu>0$ on $(a,b)$ or $\int_{[x,b]}fd\nu>0$ on $(a,b)$.
   }
\end{thm}\\

\noindent
{\bf Remarks:}

\noindent
1) In Theorem \cite{FG} and in the following, 
we say that a measure $\nu$ is $s$-affine on a segment $[a,b]$ if its density $\psi$ satisfies that 
$\psi^\gamma$ is affine on $[a,b]$, where $\gamma=\frac{s}{1-s}$.\\

\noindent
2) Notice that in Theorem \cite{FG}  it is assumed that $s\le \frac{1}{2}$. 
If $\frac{1}{2}< s\le 1$,  as follows from Theorem \cite{Bor1}, 
 the set of $s$-concave measures contains only measures whose support is one-dimensional and the Dirac measures.   
 Moreover, a quick look at the proof of Theorem \cite{FG} shows that 
the conclusions of the theorem remain valid 
except the fact that the measure $\nu$ is $s$-affine. 
It would be interesting to know if Theorem \cite{FG} may be fully extended to  $\frac{1}{2}< s\le 1$.\\

The proof of Theorem \ref{Main} splits in two steps. The first step consists in the application of 
Theorem [FG] to the reduce to the one-dimensional case and the second step is the proof of the
 one-dimensional case:\\

{\bf Step 1: Reduction to the dimension 1.}\\

Let $F$ be a Borel set in $\R^n$ and $t> 1$. Let $s\in(-\infty, 1]$ and $\mu$ be a 
$s$-concave probability such that $\mu(F_t^c)>0$. Our aim is to prove that 
$$
\mu(F^c)\ge \left(\frac{2}{t+1} \mu(F_t^c)^s + \frac{t-1}{t+1}\right)^{1/s}
\ {i.e.}\
\mu(F)\le 1-\left(\frac{2}{t+1} \mu(F_t^c)^s + \frac{t-1}{t+1}\right)^{1/s}.
$$
By density, we may assume that $\mu$ is compactly supported. 
We denote by $K$ its support which is a convex set in $\R^n$ and by $G$, 
the affine subspace generated by $K$.
Notice that in the proof of this inequality, we always may assume that 
$F\subset K$ (if we replace $F$ by $\tilde{F}:=F\cap K$, then  
$\mu(\tilde{F})=\mu(F)$ and $\tilde{F}_t\subset F_t$, hence
$\mu(\tilde{F}_t^c)\ge\mu(F_t^c)$).

From Theorem \cite{Bor1} of Borell stated in the introduction, $\mu$ is absolutely continuous 
with respect to the Lebesgue measure on $G$. Using the regularity of the measure,
we may assume that $F$ is compact in $K$. 
To satisfy the other semi-continuity hypothesis, we would need $F_t$ to be open. 
Since this is not necessarily the case, we introduce an auxiliary open set $O$ 
such that $F_t\subset O$ and $\mu(O^c)>0$.
Define $\theta\in\R$, $f:\R^n\to\R$ and $\Phi: \P(K)\to\R$ by
$$
\theta=\mu(O^c)>0,\quad f={\bf 1}_{O^c}-\theta  \quad {\rm and}\quad \Phi(\la)=\la(F),\ \forall\la\in\P(K) .
$$
Since $F$ is closed and $O$ is open, the functions $f$ and $\Phi$ are upper semi-continuous.
With these definitions, the set $P_f$ defined in the statement of the preceding theorem is
$$
P_f=\{\la\in\P(K);\ \la \ {\rm is}\ s{\rm -concave\ and}\ \la(O^c)\ge\theta\}.
$$
Since $\mu\in P_f$, if we prove that
\begin{eqnarray}\label{supphi}
\sup_{\la\in P_f}\Phi(\la)\le  1-\left(\frac{2}{t+1} \theta^s + \frac{t-1}{t+1}\right)^{1/s},
\end{eqnarray}
we will get that for any open set $O$ containing $F_t$ such that $\mu(O^c)>0$  
$$
\mu(F)\le 1-\left(\frac{2}{t+1} \mu(O^c)^s + \frac{t-1}{t+1}\right)^{1/s}.
$$
Taking the supremum on such open set $O$ and using the regularity of $\mu$, it will give the result.
From Theorem [FG], to establish inequality (\ref{supphi}) it is enough to prove it for two types of particular measure $\nu$:

- the measure $\nu$ is a Dirac measure at a point $x$ such that $f(x)\ge0$. 
It implies that ${\bf 1}_{O^c} (x)\ge\theta>0$, thus $x\notin O$, hence $x\notin F$. 
Therefore $\Phi(\delta_x)=\delta_x(F)=0$. This proves inequality (\ref{supphi}) in this case.

- the measure $\nu$ is $s$-concave on a segment $[a,b]$, such that $\int fd\nu=0$ and 
$\int_{[a,x]}fd\nu>0$ on $(a,b)$ or $\int_{[x,b]}fd\nu>0$ on $(a,b)$. 
Without loss of generality we may assume that $\int_{[x,b]}fd\nu>0$ on $(a,b)$. Hence these conditions give 
$$
\nu(O^c)=\theta \quad{\rm and}\quad \nu(O^c\cap[x,b])>\nu(O^c)\nu([x,b]),\ \forall x\in (a,b)\ .
$$
As explained at the beginning of the proof, we may assume that $F\subset[a,b]$. It is easy to see
that for a one-dimensional set $F$, its dilation $F_t$ is open in the line generated by $[a,b]$.
Hence we may choose $O=F_t$ and get rid of the auxiliary set $O$. So we have
\begin{eqnarray}\label{corel}
\nu(F_t^c)=\theta \quad{\rm and}\quad \nu(F_t^c\cap[x,b])>\nu(F_t^c)\nu([x,b]),\ \forall x\in (a,b)\ .
\end{eqnarray}
Letting $x$ tends to $b$ and using that $F_t$ is open, the second condition implies that $b\notin F_t$.
Since everything is one-dimensional, it will be more convenient to assume 
that $F\subset [a,b]\subset\R$, with $a<b$.

Let us see now why we may assume that $a\in F$. Since $F$ is closed, if $a\notin F$, then $a':=\inf F>a$. 
Let $\nu'$ be the probability, which is the restriction of $\nu$ to the interval $[a',b]$, {\em i.e.}
$\nu'=\nu_{|[a',b]}/\nu([a',b])$. Then
$$
\nu'(F^c)=\frac{\nu(F^c\cap[a',b])}{\nu([a',b])}=\frac{\nu(F^c)-\nu([a,a'])}{1-\nu([a,a'])}\le\nu(F^c)
$$
and from the second condition in (\ref{corel})
$$
\nu'(F_t^c)=\frac{\nu(F_t^c\cap[a',b])}{\nu([a',b])}>\nu(F_t^c).
$$
This ends the first step. We showed that to prove Theorem \ref{Main} 
for any $s$-concave measure $\mu$ and any Borel set $F$,  it is enough to prove it for the $s$-concave 
probabilities $\nu$ which are supported  on a segment $[a,b]\subset\R$, 
with $b\notin F_t$, $a\in F$ and $F\subset[a,b]$.
Moreover for $s\le\frac{1}{2}$, we also may assume that $\nu$ is $s$-affine.\\

{\bf Step 2: Proof in dimension 1.}\\

Let us start with a joint remark with Gu\'edon:\\
In the case where $F$ is convex, it is now easy to conclude, which 
enables us to recover the result of Gu\'edon \cite{G}.
From the convexity of $F$ and $F_t$ there exists $c<d$ such that $F=[a,c]$ and $F_t\cap [a,b]=[a,d)$ 
and we have $a<c<d<b$. Using that $d\notin F_t$ and the definition of $F_t$, for any interval 
$I$ containing $d$, we have $|I|\ge \frac{t+1}{2} |F\cap I |$. For $I=[a,d]$, this gives $d-a\ge \frac{t+1}{2} (c-a)$ and so 
$$
c\le\frac{2}{t+1}d + \left(\frac{t-1}{t+1}\right)a\quad {\rm hence}
\quad [c,b]\supset\frac{2}{t+1}[d,b] + \frac{t-1}{t+1}[a,b].
$$
Since $\nu$ is $s$-concave, we get
$$
\nu([c,b])\ge \left(\frac{2}{t+1}\nu([d,b])^s + \frac{t-1}{t+1}\nu([a,b])^s\right)^{1/s}.
$$
This ends the proof in this case since $\nu(F^c)=\nu([c,b])$, $\nu(F_t^c)=\nu([d,b])$ and $\nu([a,b])=1$. \\

The general case is more complicate. The proof of Nazarov, Sodin and Volberg \cite{NSV1}, 
to treat the log-concave ($s=0$) one-dimensional case, extends directly to the case $s\le 1$, with some suitable adaptations 
in the calculations, so we don't reproduce it here. 
But for $s\le\frac{1}{2}$, 
using that $\nu$ may be assumed $s$-affine, we can shorten the proof 
(in fact, we only use the monotonicity of the density of $\nu$). 
 
Since $F_t$ is open in $\R$, it is the countable union of disjoint intervals. By approximation, 
we may assume that there are only a finite number of them. 
Since $a\in F\subset F_t$ and $b\notin F_t$, we can write
$$
F_t\cap[a,b]=[a_0, b_0)\cup \left(\bigcup_{i= 1} ^N(a_i,b_i)\right)\quad{\rm with}\ a_i<b_i<a_{i+1},\ 0\le i\le N-1 , 
$$
where $a_0=a$. Let $F_i=F\cap (a_i,b_i)$. 
Denote by $\psi$ the density of $\nu$ with respect to the Lebesgue measure. There are two cases:\\

- If $\psi$ is non-decreasing: this is the easiest case.
Let $0\le i\le N$. Since $b_i\notin F_t$,  using the definition of $F_t$, it follows that 
for every interval $I$ containing $b_i$, we have $|I|\ge \frac{t+1}{2} \|F\cap I|$. 
Let $x\in (a_i,b_i)$, if we apply it to $I=[x,b_i]$ we get 
$$
|[x,b_i]|\ge \frac{t+1}{2} |[x,b_i]\cap F|.
$$
Hence the function $\rho:=1-\frac{t+1}{2} {\bf 1}_F$ satisfies $\int_{x}^{b_i}\rho(u)du\ge0$. Integrating by parts this gives
$$
\int_{a_i}^{b_i}\rho(u)\psi(u)du
= \psi(a_i) \int_{a_i}^{b_i}\rho(x)dx + \int_{a_i}^{b_i}\left(\int_{x}^{b_i}\rho(u)du\right)\psi'(x)dx\ge 0.
$$
Hence $\nu\big((a_i, b_i)\big)\ge \frac{t+1}{2} \nu(F_i)$ and since $F_t=\cup (a_i,b_i)$ and $F=\cup F_i$, it follows that
$\nu(F_t)\ge \frac{t+1}{2} \nu(F)$. Therefore, using the comparison between the $s$-mean (with $s\le 1$) 
and the arithmetic mean, we conclude that
$$
\nu(F^c)\ge \frac{2}{t+1} \nu(F_t^c) + \frac{t-1}{t+1}\ge \left(\frac{2}{t+1} \nu(F_t^c)^s + \frac{t-1}{t+1}\right)^{1/s}.
$$

- If $\psi$ is non-increasing: 
We first prove that, for each $0\le i\le N$
\begin{eqnarray}\label{fi}
\nu(F_i^c)\ge \left(\frac{2}{t+1} \nu((a_i,b_i)^c)^s + \frac{t-1}{t+1}\right)^{1/s}.
\end{eqnarray}
For $i\ge 1$, we have $a_i\notin F_t$ and it is similar as the previous case. Indeed, for every $x\in(a_i,b_i)$,
$|[a_i,x]|\ge \frac{t+1}{2} |[a_i,x]\cap F|$ and an integration by parts gives that $\nu\big((a_i, b_i)\big)\ge \frac{t+1}{2} \nu(F_i)$.
From the comparison of the means, inequality (\ref{fi}) follows.\\
For $i=0$, we have $a_0=a\in F$. We define $F_0'=[a_0,c_0]$, where $c_0$ is chosen such that
$|F_0'|=|F_0|$. Since $\psi$ is non-increasing, we have $\nu(F_0')\ge\nu(F_0)$ and since $b_0\notin F_t$,
$$
|[a_0,b_0]|\ge \frac{t+1}{2} |[a_0,b_0]\cap F|=\frac{t+1}{2} |F_0|=\frac{t+1}{2} |F_0'|=\frac{t+1}{2} |[a_0,c_0]|.
$$
Hence $b_0-a_0\ge \frac{t+1}{2} (c_0-a_0)$. As in the joint remark with Gu\'edon given before, 
we get that
$$
\nu([c_0,b])\ge \left(\frac{2}{t+1}\nu([b_0,b])^s +\frac{t-1}{t+1}\nu([a_0,b])^s\right)^{1/s}.
$$
Therefore we get inequality (\ref{fi}) for $i=0$:
$$
\nu(F_0^c)\ge\nu(F_0'^c)=\nu([c_0,b])\ge\left(\frac{2}{t+1}\nu((a_0,b_0)^c)^s + \frac{t-1}{t+1}\right)^{1/s}.
$$
The inequality (\ref{fi}) may be written $\nu(F_i)\le \f(\nu\big(a_i,b_i)\big)$, for $0\le i\le N$, 
where $\f:[0,1]\to [0,1]$ is defined by 
$$
\f(x)=1-\left(\frac{2}{t+1}(1-x)^s + \frac{t-1}{t+1}\right)^{1/s}.
$$
From Minkowski inequality for the $s$-mean, with $s\le 1$, the function $\f$ is convex on $[0,1]$.
Denote $\la_i=\nu(\big(a_i,b_i)\big)/\nu(F_t)$. Using that $\f(0)=0$ and the convexity of $\f$ we get 
$$
\nu(F_i)\le \f(\nu\big(a_i,b_i)\big)=\f\big(\la_i\nu(F_t)\big)\le\la_i\f\big(\nu(F_t)\big).
$$
Summing on $i$ and using that $\sum_{i=1}^N\la_i=1$, we conclude that 
$$
\nu(F)\le\f\big(\nu(F_t)\big).
$$
This is the result.
\hfill\qed\\

\noindent
{\bf Acknowledgments:} The author thanks Olivier Gu\'edon and Jean Saint Raymond for useful discussions.

\end{document}